\documentclass[12pt,a4paper]{article}
\usepackage{amsmath}
\usepackage{amssymb}
\usepackage{bm}
\usepackage[all]{xy}

\usepackage{amsfonts,amsthm,euscript}

\theoremstyle{remark}
\newtheorem*{rem*}{Remark}

\def\genfd{{\bf k}} 
\def\id{{\rm id}}

\newcounter{point}
\newcommand{\ppt}{\vskip .01in\addtocounter{point}{1}{\bf \arabic{point} }}
\newcommand{\ppta}{{\bf \arabic{point}a }}
\newcommand{\pptb}{{\bf \arabic{point}b }}
\newcommand{\pptc}{{\bf \arabic{point}c }}

\def\lnamedef#1{\expandafter\edef\csname lpt#1 \endcsname}
\def\lnameuse#1{\expandafter\csname lpt#1 \endcsname}
\newcommand\ldef[1]{\lnamedef{#1}{\arabic{point}}}
\newcommand\luse[1]{{\bf\lnameuse{#1}}}
\newcommand{\nodo}[1]{}

\begin{document}

\title{Distributive laws for actions of monoidal categories}

\author{
Zoran \v{S}koda
}

\date{{}}

\maketitle
\begin{abstract}
\noindent Given a monoidal category $\cal C$, 
an ordinary category $\cal M$, and a monad $\bf T$ in $\cal M$,  
the lifts in a strict sense of a fixed action of 
${\cal C}$ on ${\cal M}$ to an action of 
$\cal C$ on the Eilenberg-Moore category ${\cal M}^{\bf T}$ 
of $\bf T$-modules are in a bijective correspondence 
with certain families of natural transformations,
indexed by the objects in the monoidal category $\cal C$.   
These families are analogues of
distributive laws between monads. 
\end{abstract}

\ppt Actions of monoidal categories on other ordinary categories
appear in diverse setups. Anyway, 
apart from explaining the main result of this paper,
we will limit ourselves just to our original motivation,
and other applications will be left to the reader's
imagination and future.

\ppt \ldef{standard}
Given some group-like (symmetry) object, e.g.
a Hopf algebra $H$ (e.g $H= {\cal O}(G)$ for an algebraic group $G$),
and an auxiliary category $\cal M$
viewed as a category ${\frak Qcoh}_X$ 
of (i.e. ``quasicoherent'') sheaves on a (noncommutative) space $X$,
what is the appropriate notion of the action of $G$ on $X$
and how do we represent it practically in such a situation ?
Standard observations: 

\ppta {\it Affine case}. If $\cal M$ is the category of
(say left) modules over an algebra $A$, 
then in {\sc Drinfeld}'s quantum group philosophy 
one replaces group action by Hopf (co)actions 
of some Hopf algebra on $A$. 

\pptb {\it (Co)monads from (co)actions.} 1) {\it Affine case}.
Every Hopf coaction of bialgebra $H$ on an algebra $A$
induces a comonad, say $\bf T$, on $A-{\rm Mod}$.

2) {\it Globalizing}. In some cases one may globalize 
this picture by introducing a global (co)monad $\bf T'$ and having
some device which will locally compare/reduce its action to the action
induced from  ``affine Hopf-coaction'' (as in a)).
In the commutative situation, and if $H = {\cal O}(G)$ is 
(a Hopf algebra of functions on a group), one can localize the 
space with action only to the invariant open sets, otherwise
there is no induced action. Although we have used 
generalizations of such methods to noncommutative setup
elsewhere, {\em a priori} there is no sufficient supply 
of affine open sets which are action-invariant. 

\pptc One occasionally replaces the Hopf algebra 
(or more general symmetry objects) by a monoidal category 
$\tilde{\cal C} = ({\cal C},\otimes)$
of its (co)modules. Say, in the case of modules, 
the Hopf algebra itself is a distinguished comonoid in $\tilde{\cal C}$.
Its image under the (action) functor ${\cal L}:{\cal C}\to{\rm End}\,{\cal M}$
is the comonad $\bf T$ we mentioned before. 
The category ${\cal M}^{\bf T}$ of modules over this monad 
is a version of the category of $G$-equivariant sheaves on $X$ 
in this setup.

\ppt The crucial point missing from \luse{standard}{\bf a-c}, 
and from the bulk of recorded ``noncommutative'' literature is 
failing to adhere to the old {\sc Grothendieck}'s advice
that geometrical study of morphisms (hence actions in particular)
should be carried in {\it relative setup}. 
In particular, {\it group schemes should live over some fixed scheme $V$ 
and the objects they act on should live over that same scheme $V$}. 
Similarily, Hopf algebras
are defined in some (symmetric, but generalizations apply) 
monoidal category ${\cal V}$. The space $X = ``{\rm Spec}\,{\cal M}''$ 
should live over some ``${\rm Spec}\,{\cal V}$'' 
what means that we have direct image
(forgetful) functor from ${\cal M}$ to ${\cal V}$, 
satisfying some properties (e.g. having an inverse image functor)
needed to call $X$ a space or a scheme
(e.g. a noncommutative scheme over $\cal V$ (\cite{Ros:NcSch})).
The action has to respect the forgetful functor
in the sense that if one applies the action and then forgets
down to ${\cal V}$ then we (should) obtain exactly the natural
action of the Hopf algebra (i.e. of its category of (co)modules) on
${\cal V}$. In this paper, we show a general nonsense result
on what additional data are necessary to lift known actions
of monoidal category on ${\cal V}$ 
(e.g. the natural action of $H-{\rm Comod}$ on
${\rm Vec}_\genfd$ where $H$ is a bialgebra over a field $\genfd$)
to new actions on a simple, but important class of categories over $\cal V$:
namely the categories of modules over monads in $\cal V$. 
Fundamental importance of such categories 
in geometric context (as local models) 
has been enlightened in~\cite{Ros:NcSch}. 

Correct choice of $\cal V$ 
enables unified correct treatment of diverse flavours of symmetries, 
e.g. relative group schemes, Hopf algebras in braided categories,
(finite) group algebras in monoidal categories different 
than (super)vector spaces, accounting for (multi)graded rings etc.
Useful sources of actions in quantum vector bundle theory, e.g. 
``entwining structures'' (\cite{BrzMajid:ent,Bohm:ent}) 
naturally appear in our perspective. 
Natural generalizations work, primarily for {\it bialgebroids}.
The relative viewpoint is also inherent in replacing Hopf algebras
by tensor categories in view of Tannakian theory which does require
a fibre functor. For sensible consideration of actions 
one should not detach a monoidal category from its origin.

Some relevant references for the geometrical motivation for this work
are~\cite{BrzMajid:ent, LuntsRosMP,LuntsSkoda:desc, 
Ros:NcSS, Skoda:coh-states}. 
\vskip .14in

\ppt A {\bf monoidal category} is a 6-tuple 
$\tilde{\cal C} = ({\cal C},\otimes, {\bf 1}, a, r, l)$ where

$\otimes : {\cal C}\times {\cal C}\to {\cal C}$ is a bifunctor 
(monoidal or tensor product),

${\bf 1} \in {\rm Ob}\,{\cal C}$ is a distinguished object (unit object),

$a : \_ \otimes (\_\otimes\_) \Rightarrow (\_\otimes\_)\otimes\_$
is a trinatural equivalence of trifunctors (associativity constraint);
(i.e. a family of isomorphisms $\{a_{X,Y,Z} : X\otimes (Y\otimes Z) 
\to (X\otimes Y)\otimes Z \}_{X,Y,Z\in {\rm Ob}\,{\cal C}}$,
natural in $X,Y,Z$),

$\rho : {\rm Id}_{\cal C} \to {\rm Id}_{\cal C}\otimes {\bf 1}$
(right unit coherence) and 
$\lambda : {\rm Id}_{\cal C} \to {\bf 1} \otimes {\rm Id}_{\cal C}$
(left unit coherence) are
natural equivalences of functors; the compositions\newline
$X\otimes(Y\otimes(Z\otimes W))\stackrel{a_{X,Y,Z\otimes W}}\longrightarrow
(X\otimes Y)\otimes(Z\otimes W)\stackrel{a_{X\otimes Y,Z,W}}\longrightarrow
((X\otimes Y)\otimes Z)\otimes W$
and $X\otimes(Y\otimes(Z\otimes W))
\stackrel{X\otimes a_{Y,Z,W}}\longrightarrow
X\otimes ((Y\otimes Z)\otimes W)\stackrel{a_{X,Y\otimes Z,W}}
\longrightarrow(X\otimes (Y\otimes Z))\otimes W
\stackrel{a_{X,Y,Z}\otimes W}
\longrightarrow ((X\otimes Y)\otimes Z)\otimes W$
should agree $\forall X,Y,Z,W\in{\rm Ob}\,{\cal C}$ and
the triangle coherence identity 
$a_{X,{\bf 1}, Y} \circ (\rho_X \otimes Y) = X\otimes \lambda_Y$ holds.
The monoidal category is {\bf strict} if
all coherence isomorphisms 
$a_{X,Y,Z}$, $\rho_X$, $\lambda_X$ are identity maps
and hence may be omitted from the data. 

\ppta A {\bf monoidal functor} from a monoidal category
$({\cal A},\otimes,{\bf 1},a,\rho,\lambda)$
to $({\cal B},\otimes',{\bf 1}',a',\rho',\lambda')$
is a triple $(F,\chi,\xi)$ where $F:{\cal A}\to {\cal B}$ is a functor
and $\chi:F(\_)\otimes' F(\_)\to F(\_\otimes\_)$ and
$\xi : F(\_)\otimes'{\bf 1}'\to F(\_\otimes {\bf 1})$ 
natural transformations satisfying some natural coherence relations
~(\cite{MacLane, Bunge:relfunct}).

\ppt The category ${\rm End}\,{\cal A}$ 
of endofunctors in a given category ${\cal A}$
is a strict monoidal category
with respect to the tensor product 
given on objects (endofunctors) as the composition,
and on morphisms (natural endotransformations) 
as the Godement's ``vertical product'' 
$\star : (F:f\Rightarrow f', G:g\Rightarrow g')\mapsto 
G\star F = : g\circ f\Rightarrow g'\circ f'$,
where $(G\star F)_M := g'(F_M)\circ G_{f(M)} 
= G_{f'(M)}\circ g(F_M)$. 
A semigroup ${\bf T} = (T,\mu)$ in $({\rm End}\,{\cal A},\circ,{\rm Id})$ 
is called a {\bf nonunital monad}.
In other words, $T$ is an endofunctor in ${\cal A}$
and $\mu : T\circ T\to T$ is a natural transformation of functors
satisfying the associativity 
$({\rm Id}\star T)\circ T = (T\star{\rm Id})\circ T$.
If $(T,\mu)$ is a monoid in ${\rm End}\,{\cal A}$,
i.e. in addition $\exists$ (automatically unique) natural
transformation $\eta : {\rm Id}\Rightarrow T$ satisfying
$\mu\circ\eta_T = {\rm Id} = \mu\circ T(\eta)$, then we say that
$(T,\mu,\eta)$ is a {\bf monad} (\cite{MacLane}). 
A {\bf comonad} is a monad in ${\cal A}^{\rm op}$.

\ppt A right {\bf action} of a monoidal category is a bifunctor
\[ \lozenge : {\cal M}\times {\cal C}\to {\cal M}, \,\,\,\,\,\,\,\,\,
(M,Q)\mapsto M \lozenge Q, \]
and a family of morphisms $\Psi_M^{X,Y} : M \lozenge (X\otimes Y)
\mapsto (M\lozenge X)\lozenge Y$ in ${\cal M}$ natural
in $M, X$ and $Y$ and such that the ``action associativity'' pentagon
$$\xymatrix{
M \lozenge (X \otimes (Y \otimes Z)) 
\ar[rr]^{\Psi_M^{X,Y\otimes Z}} 
\ar[d]_{M\lozenge\, a^{-1}_{X,Y,Z}} &&  
(M\lozenge X) \lozenge (Y \otimes Z)\ar[dd]^{\Psi_{M\lozenge X}^{Y,Z}} \\
M \lozenge ((X \otimes Y)\otimes Z)\ar[d]_{\Psi_M^{(X\otimes Y),Z}}  &&\\
(M \lozenge (X\otimes Y)) \lozenge Z \ar[rr]^{\Psi_M^{X,Y}\lozenge Z}  &&
((M\lozenge X)\lozenge Y)\lozenge Z 
}$$
commutes. For applications we have in mind, we do {\it not} 
require that the unit object $\bf 1$ acts strictly trivially,
i.e. $M \lozenge {\bf 1}\neq M$ in general,
but rather we demand a natural equivalence of functors
$u = u^\lozenge: {\rm Id} \to {\rm Id}\lozenge {\bf 1}$, 
compatible with coherence isomorphisms in $\tilde{\cal C}$.
More precisely, the diagram
$$\xymatrix{
(M\lozenge{\bf 1})\lozenge Q &&
M\lozenge({\bf 1}\otimes Q)\ar[ll]_{\Psi_M^{{\bf 1},Q}}
\ar[dl]^{M\lozenge l_Q}
\\&M\lozenge Q\ar[ld]_{u_{M\lozenge Q}}\ar[lu]^{(u_M)\lozenge Q}&\\
(M\lozenge Q)\lozenge{\bf 1}&& 
M\lozenge(Q\otimes{\bf 1})\ar[ll]_{\Psi_M^{Q,{\bf 1}}}
\ar[ul]_{M\lozenge r_Q}
}$$
commutes for all $M\in {\rm Ob}\,{\cal M}$ and all $Q\in {\rm Ob}\,{\cal C}$.

\ppta A right {\bf action} of a monoidal category is 
$({\cal C},\otimes, a, r, l)$ on a 
(non-monoidal) category ${\cal M}$ may be also given by a contravariant
monoidal functor from ${\cal C}$ to ${\rm End}\,{\cal M}$.

Given an action $(\lozenge, \Psi, u)$,
the associated coherent monoidal functor 
${\cal L} : {\cal C}\to {\rm End}\,{\cal M}$, 
is then given by ${\cal L}(C)(M) :=  {\cal L}_C(M) := M\lozenge C$ (the rest
of the coherence structure left to the reader). 

\pptb {\it Example.} Every monoidal category $\tilde{\cal C} = 
({\cal C}, \otimes, {\bf 1}, a,r,l)$ acts on its underlying category $\cal C$
from the right by the action $(\lozenge,\Psi, u) = (\otimes, a,r^{-1})$.

\vskip .01in
\ppt A {\bf distributive law} from a monoidal action $\lozenge$ of $\cal C$
on $\cal M$ to a monad ${\bf T} = (T,\mu,\eta)$ is a binatural
transformation $l : T(\_\lozenge\_)\Rightarrow(T\_)\lozenge\_$
of bifunctors ${\cal M} \times {\cal C} \to {\cal M}$,
satisfying the list of axioms (D1-4) below.
Denote by $(l^Q)_M = l^Q_M := l_{M,Q} : T(M\lozenge Q) \to TM\lozenge Q$ the
morphisms in $\cal M$ forming $l$. 
Clearly, each $l^Q : T(\_\lozenge Q)\Rightarrow T(\_)\lozenge Q$ 
is a natural transformation, 
hence $l$ may be viewed as a family $\{l^Q\}_{Q\in {\rm Ob}\,{\cal C}}$.
The diagrams (D1-4) are required to commute:
\begin{equation}\label{D1}
\xymatrix{
TT(M \lozenge Q)\ar[r]^{T(l^Q_M)}
\ar[d]_{\mu_{M\lozenge Q}} &
T (T(M)\lozenge Q) \ar[r]^{l^Q_{TM}} &
  TT(M)\lozenge Q
\ar[d]^{\mu_M\lozenge Q}\\
T(M\lozenge Q)\ar[rr]^{l^Q_M}&& T(M)\lozenge Q.
}
\tag{D1}
\end{equation}
\begin{equation}\label{D2}
\xymatrix{
T(M\lozenge (Q\otimes Q')) 
\ar[rr]^{l^{Q\otimes Q'}_M}\ar[d]^{T(\Psi _M^{Q,Q'})}
&& TM\lozenge (Q\otimes Q')\ar[d]^{\Psi_M^{Q,Q'}}\\
T((M\lozenge Q)\lozenge Q') \ar[r]^{l^{Q'}_{M\otimes Q}}
& T(M\lozenge Q)\lozenge Q' \ar[r]^{l^Q_M\lozenge Q'}
& (TM\lozenge Q)\lozenge Q'
}\tag{D2}\end{equation}
\begin{equation}\label{D3}\xymatrix{
& M\lozenge Q\ar[ld]_{\eta_{M\lozenge Q}}\ar[rd]^{\eta_M\lozenge Q} &\\
T(M\lozenge Q)\ar[rr]^{l^Q_M} && T(M)\lozenge Q 
}\tag{D3}\end{equation}
\begin{equation}\label{D4}\xymatrix{ 
& T(M)\ar[ld]_{T(u_{M})}\ar[rd]^{u_{TM}} &\\
T(M\lozenge 1)\ar[rr]^{l^Q_M} && TM\lozenge 1 
}\tag{D4}\end{equation}

\ppt Given two $\cal C$-categories $({\cal M},\lozenge,\Psi, u)$,
$(\tilde{\cal M},\tilde\lozenge,\tilde\Psi,\tilde{u})$,
a pair $(F,\zeta)$ where $F : \tilde{\cal M}\to \cal M$ is a functor, 
and $\zeta : F(\_)\tilde\lozenge\_ \Rightarrow F(\_\lozenge\_)$
a binatural equivalence of bifunctors 
$\tilde{\cal M}\times{\cal C}\to {\cal M}$,
is called a $\cal C$-{\bf functor} if the following diagrams commute
$$\xymatrix{
F(M) \lozenge {\bf 1} \ar[rr]^{\zeta_{M,{\bf 1}}}\ar[rd]_{u_{F(M)}} 
&& F(M\tilde\lozenge {\bf 1})\ar[ld]^{F(\tilde{u}_M)} \\
& F(M) &
}$$
$$\xymatrix{
F(M)\lozenge(Q\otimes Q')
\ar[r]_{\Psi^{Q,Q'}_{F(M)}}\ar[d]_{\zeta_{M,Q\otimes Q'}}&
(F(M)\lozenge Q)\lozenge Q'\ar[r]_{\zeta_{M,Q}\lozenge Q'}&
F(M\tilde\lozenge Q)\lozenge Q'\ar[d]^{\zeta_{M\tilde\lozenge Q, Q'}}\\
F(M\tilde\lozenge(Q\otimes Q'))\ar[rr]^{F(\tilde\Psi_M^{Q,Q'})}&&
F((M\tilde\lozenge Q)\tilde\lozenge Q')
}$$
We say that $\cal C$-functor $(F,\zeta)$ is {\bf strict} if 
$\zeta_{M,Q} = {\rm Id}_{F(M\tilde\lozenge Q)}$ for all $M$ and $Q$.
Then we naturally omit $\zeta$ from the notation. 

\ppt If $U : \tilde{\cal M}\to{\cal M}$ is some 
(usually ``forgetable'') functor,
then a $\cal C$-category structure $(\tilde{\cal M},\tilde\Psi,\tilde{u})$
is called a {\bf strict lift} of a ${\cal C}$-category
structure $({\cal M},\Psi,u)$ to $\tilde{\cal M}$ along $U$ 
if $U$ is a strict ${\cal C}$-functor
from $(\tilde{\cal M},\tilde\Psi,\tilde{u})$ to $({\cal M},\Psi,u)$.
In other words, $U(M\tilde\lozenge Q) = U(M)\lozenge Q$,
$u_{U(M)} = U(\tilde{u}_M)$ 
and $U(\tilde\Psi_M^{Q,Q'}) = \Psi_{U M}^{Q,Q'}$.

\ppt In this article we are primarily interested in the case when 
$\tilde{\cal M} = {\cal M}^{\bf T}$ is the {\bf Eilenberg-Moore category}
of a monad ${\bf T}$ in ${\cal M}$. Its objects are $\bf T$-{\bf modules},
that is pairs of the form $(M,\nu)$, 
where $\nu : TM \to M$ is a $\bf T$-{\bf action}, 
i.e. a morphism satisfying $\nu \circ T(\nu) = \nu\circ\mu_M$
and the unitality $\nu\circ\eta_M = \id_M$; 
the morphisms $f : (M,\nu)\to (N,\xi)$ are those morphisms $f : M \to N$  
which intertwine the action in the sense $\xi \circ T(f) = f \circ \nu$.
There is a forgetful functor $U = U^{\bf T} : {\cal M}^{\bf T} \to {\cal M}$
given by $U^{\bf T} : (M,\nu)\mapsto M$. This forgetful functor
will be in place of functor $U : \tilde{\cal M}\to {\cal M}$ 
in previous discussion.
$U^{\bf T}$ has a left adjoint $F^{\bf T} : M\mapsto (TM,\mu_M)$ 
whose image is so-called Kleisli category
of a monad ${\bf T}$ and its objects are called free $\bf T$-modules. 
Note that $U^{\bf T}F^{\bf T} = T$. 

\ppt \ldef{prellift} Suppose some action 
$(\tilde\lozenge,\tilde\Psi,\tilde{u})$
strictly lifts $(\lozenge, \Psi, u)$ to ${\cal M}^{\bf T}$. 
Then we may write 
$(M,\nu)\tilde\lozenge Q = (M\lozenge Q,\widehat{(M,\nu)}^Q)$ where
$\widehat{(M,\nu)}^Q : T(M\lozenge Q)\to (M\lozenge Q)$ is a
(canonical) $\bf T$-action (depending on $M,\nu,Q$). Rule
$(M,\nu)\mapsto \widehat{(M,\nu)}^Q$ defines a
natural transformation $\hat{()}^Q : T(U(\_)\lozenge Q)\to U(\_)\lozenge Q$
of functors from ${\cal M}^{\bf T}$ to ${\cal M}$.
We will sometimes abuse the notation abbreviating $\hat{\nu}^Q$ for
$\widehat{(M,\nu)}^Q$. The counit of the adjunction 
$\epsilon : F^{\bf T}U^{\bf T}\Rightarrow {\rm Id}_{{\cal M}^{\bf T}}$
is given by $\epsilon_{(M,\nu)} = \nu : F^{\bf T}M = (TM,\mu_M)\to (M,\nu)$.
Then $U^{\bf T}(\epsilon_{(M,\nu)}) = \nu : TM\to M$.
Thus $\widehat{(M,\nu)}^Q = U^{\bf T}(\epsilon_{(M,\nu)\tilde\lozenge Q})$, 
i.e. $\hat{()}^Q = U^{\bf T}(\epsilon_{\_\tilde\lozenge Q})$.
We may also write
$$\hat{()}^Q F^{\bf T} = U^{\bf T}(\epsilon_{F^{\bf T}(\_)\tilde\lozenge Q}) : 
M\mapsto \left(\widehat{(TM,\mu_M)}^Q 
: T(TM\lozenge Q)\to TM\lozenge Q\right)$$ 
what is a natural transformation of functors ${\cal M}\to {\cal M}$.

\ppt {\bf Theorem.} {\it
Distributive laws from an action $(\lozenge, \Psi, u)$ of 
a monoidal category $\tilde{\cal C} = ({\cal C},\otimes,{\bf 1},a,r,l)$ 
on a category ${\cal M}$ 
to a monad ${\bf T} = (T,\mu,\eta)$ in ${\cal M}$ are 
in a natural bijective correspondence
$l \mapsto (\lozenge_l,\Psi_l,u_l)$
with the monoidal actions of ${\cal C}$ on ${\cal M}^{\bf T}$ 
strictly lifting the action $(\lozenge, \Psi, u)$ of ${\cal C}$ on ${\cal M}$
along the forgetful functor $U = U^{\bf T}:{\cal M}^{\bf T}\to {\cal M}$,
$U^{\bf T} : (M,\nu) \mapsto M$.
} 

{\it Proof.} I. From distributive laws to lifted actions:
$l \mapsto  (\lozenge_l,\Psi_l, u_l)$.

Given a distributive law $l$, 
the action of $\cal C$ on ${\cal M}$ is lifted to
an action of ${\cal C}$ on ${\cal M}^{\bf T}$
which is of the form $(M,\nu)\lozenge_l Q := (M\lozenge Q,\nu^Q)$,
where $\nu^Q$ is the composition
$$\xymatrix{
T(M\lozenge Q)\ar[rr]^{l^Q_M}&& TM \lozenge Q 
\ar[rr]^{\nu \lozenge Q}&& M \lozenge Q.
}$$
Actually, to live up to our claims we extend these formulas 
to a bifunctor $\lozenge_l$, using the more obvious formulas on morphisms
\[\left.\begin{array}{l}
f \lozenge_l Q := f \lozenge Q\\
f \lozenge_l g := f \lozenge g
\end{array}\right\rbrace\begin{array}{l}
\forall f : (M,\nu) \to (N,\xi), \\
\forall g:Q \to Q'.\end{array}
\]
Pair $(M\lozenge Q,\nu^Q_M)$ is indeed a $\bf T$-module by the 
commutativity of diagram
$$\xymatrix{
TT(M\lozenge Q) \ar[r]^{T(l^Q_M)} \ar[dd]_{\mu_{M\lozenge Q}}&
T(TM \lozenge Q) \ar[r]^{T(\nu\lozenge Q)}\ar[d]^{l^Q_{TM}} &
T(M \lozenge Q) \ar[d]^{l^Q_M}\\
& TTM \lozenge Q \ar[r]^{T\nu\lozenge Q}\ar[d]^{\mu_M\lozenge Q} &
TM \lozenge Q\ar[d]^{\nu\lozenge Q} \\
T(M \lozenge Q) \ar[r]^{l^Q_M} & TM\lozenge Q\ar[r]^{\nu\lozenge Q} &
M\lozenge Q.
}$$
The pentagon on the left is (D1), 
the upper-right square expresses the naturality of $l^Q$ and
the lower-right square is obtained from the action axiom for $\nu$
by applying $\lozenge Q$. 
The action $\nu^Q : T(M\lozenge Q)\to M\lozenge Q$
is unital as 
$\nu^Q \circ \eta_{M\lozenge Q} = (\nu\lozenge Q) \circ l^Q_M 
\circ \eta_{M\lozenge Q} = (\nu\lozenge Q) \circ (\eta_M\lozenge Q)
= (\nu\circ \eta_M)\lozenge Q = \id_M \lozenge Q = \id_{M\lozenge Q}$.

If $f : (M,\nu)\to (N,\xi)$ then
$f \lozenge Q$ in ${\cal M}$ is indeed a $\bf T$-module
map from $(M\lozenge Q,\nu^Q)$ to $(N\lozenge Q,\xi^Q)$, i.e.
$(f\lozenge Q)\nu = T(f\lozenge Q)\xi^Q$ 
by the commutativity of the following diagram:
$$\xymatrix{
T(M\lozenge Q) \ar@/^/[rr]^{\nu^Q}
\ar[r]_{l^Q_M} \ar[d]_{T(f\lozenge Q)} & 
TM \lozenge Q \ar[r]_{\nu\lozenge Q} \ar[d]_{Tf\lozenge Q} &
M \lozenge Q \ar[d]_{f\lozenge Q}\\
T(N\lozenge Q) \ar@/_/[rr]_{\xi^Q}
\ar[r]^{l^Q_N} & 
TN \lozenge Q \ar[r]^{\xi\lozenge Q}  &
N \lozenge Q.
}$$
Here the left square is commutative by the naturality of $l^Q$ 
and the right square by the functoriality 
of $\lozenge$ in the first variable.

$u_l : {\rm Id}\lozenge_l {\bf 1}\Rightarrow {\rm Id}$ is defined
to be given simply by $(u_l)_{(M,\nu)} := u_M$. It is indeed a map
of $\bf T$-modules $(u_l)_{(M,\nu)} : (M,\nu)\lozenge_l {\bf 1} = 
(M\lozenge {\bf 1},\nu^{\bf 1})\to (M,\nu)$ by axiom (D4) and 
naturality of $u$:
$$\xymatrix{
TM \ar[rr]^{T(u_M)}\ar[dd]^{\nu}\ar[rd]^{u_{TM}} && 
T(M\lozenge 1) \ar[dl]^{l^{\bf 1}_M}\ar[dd]^{\nu^{\bf 1}} \\
&TM\lozenge 1\ar[dr]^{\nu\lozenge {\bf 1}}&\\
M\ar[rr]^{u_M}&&M\lozenge {\bf 1}.
}$$

The required associativity maps 
of $\bf T$-modules $(\Phi_l)^{Q,Q'}_{(M,\nu)} : 
((M,\nu) \lozenge_l (Q\otimes Q'),\nu^{Q\otimes Q'})
\to (((M,\nu)\lozenge_l Q)\lozenge_l Q',(\nu^Q)^{Q'})$ 
are defined to be identical to the maps of underlying objects in ${\cal M}$,
$\Phi^{Q,Q'}_M : M \lozenge (Q\otimes Q')\to (M\lozenge Q)\lozenge Q'$,
but one has to check that they are in fact maps of the $\bf T$-modules:
$$\xymatrix{
T(M\lozenge (Q \otimes Q')) 
\ar@/^/[rrr]^{\nu^{Q\otimes Q'}}
\ar[rr]_{l^{Q\otimes Q'}_M} 
\ar[d]^{T\Psi^{Q,Q'}_M} &&
TM \lozenge (Q \otimes Q') \ar[r]_{\nu\lozenge (Q\otimes Q')}
\ar[d]^{\Psi^{Q,Q'}_{TM}} &
M \lozenge (Q \otimes Q') \ar[d]^{\Psi^{Q,Q'}_M} \\
T((M\lozenge Q)\lozenge Q')\ar[r]^{l^Q_{M\lozenge Q}}
\ar@/_/[rrr]_{(\nu^Q)^{Q'}}
& T(M\lozenge Q)\lozenge Q' \ar[r]^{l^Q_M\lozenge Q'}&
(TM\lozenge Q)\lozenge Q' \ar[r]^{(\nu\lozenge Q)\lozenge Q'}
& (M \lozenge Q)\lozenge Q'.
}$$
The left square is commutative by the definition of the distributive
laws and the right by the naturality of $\Psi^{Q,Q'}$. 
Hence the commutativity of the external (round) square
follows, expressing the fact that
$\Psi^{Q,Q'}_M$ is a $\bf T$-module map. 
Consequently, 
setting $\tilde\Psi^{Q,Q'}_{(M,\nu)} := \Psi^{Q,Q'}_M$ is meaningful;
and the action pentagon for $\tilde\Psi$ is automatically commutative
in category $\tilde{\cal M}$; the same for compatibilities with $u_l$.

II. From lifted actions to distributive laws: 
$(\tilde\lozenge, \tilde{\Psi}, \tilde{u})\mapsto l$.

We will repeatedly use the notation and properties from~\luse{prellift}.

The action axiom ensures that any $\bf T$-action $\nu$ is a map
of $\bf T$-modules $\nu : (TM,\mu_M)\to (M,\nu)$. 
As $\tilde\lozenge$ is a bifunctor, the map
$\nu \tilde\lozenge Q$ is also a map of $\bf T$-modules,
$$
\nu \tilde\lozenge Q : (TM,\mu_M)\tilde\lozenge Q \to (M,\nu)\tilde\lozenge Q.
$$
That means that the following diagram in $\cal M$ commutes:
$$\xymatrix{
T(TM \lozenge Q) \ar[r]^{\widehat{\mu_M}^Q}\ar[d]_{T(\nu\lozenge Q)} & 
TM\lozenge Q \ar[d]^{\nu\lozenge Q}\\
T(M\lozenge Q) \ar[r]^{\hat\nu^Q} & M\lozenge Q,
}$$
where we also used that $U^T (\nu\tilde\lozenge Q) = (\nu\lozenge Q) U^T$.
This implies that 
$$\begin{array}{lcl}
(\nu\lozenge Q)\circ \widehat{\mu_M}^Q \circ T(\eta_M\lozenge Q) & = &
\hat\nu^Q \circ T(\nu\lozenge Q) \circ T(\eta_M\lozenge Q) \\
& = &\hat\nu^Q \circ T((\nu\circ\eta_M)\lozenge Q) \\
& = &\hat\nu^Q \circ T(\id_M \lozenge Q)\\
& = &\hat\nu^Q \circ \id_{T(M\lozenge Q)}\\
& = &\hat\nu^Q.
\end{array}$$ 
Therefore if we set 
$$\xymatrix{
 T(M\lozenge Q)\ar@/_/[rrrrr]_{l^Q_M}\ar[rr]^{T(\eta_M\lozenge Q)}&&
 T(TM\lozenge Q)\ar[rrr]^{\widehat{\mu_M}^Q = 
U(\epsilon_{FM\tilde\lozenge Q})} 
&&& TM\lozenge Q.
}$$
that will define a candidate for a distributive law 
$l = \{l^Q\}_{Q\in {\cal C}}$ for which
$\hat\nu^Q = (\nu\lozenge Q) \circ l^Q_M \equiv \nu^Q$ for all 
$\bf T$-modules $(M,\nu)$ (cf. the notation $\nu^Q$ from part I). 
{\bf This means, in effect, that
the lift $\lozenge_l Q$ agrees with $\tilde{\lozenge}Q$ (*)}.
We have to check that this family is really
a distributive law. The fact that $l^Q_M$ form a natural transformation 
$l^Q$ follows as for every $f : M \to N$ in $\cal M$ the diagram
$$\xymatrix{
T(M\lozenge Q)\ar[rr]^{T(\eta_M\lozenge Q)}\ar[d]_{T(f\lozenge Q)}&&
T(TM\lozenge Q)\ar[r]^{\widehat{\mu_M}^Q}\ar[d]_{T(Tf\lozenge Q)}&
TM\lozenge Q\ar[d]_{Tf}\\
T(N\lozenge Q)\ar[rr]^{T(\eta_N\lozenge Q)}&&
T(TN\lozenge Q)\ar[r]^{\widehat{\mu_N}^Q} & TN\lozenge Q
}$$
commutes, where we identified $\tilde\mu_M$ with the corresponding map
in $\cal M$.

Axiom (D3) reads $\eta_M\lozenge Q = l^Q_M \eta_{M\lozenge Q}
= U(\epsilon_{FM\tilde\lozenge Q}) T(\eta_M\lozenge Q)\eta_{M\lozenge Q}$
hence it follows by the commutativity of the pentagon
$$\xymatrix{
M\lozenge Q\ar[rr]^{\eta_M\lozenge Q}\ar[d]_{\eta_{M\lozenge Q}} &&
TM\lozenge Q\ar[d]^{\eta_{TM\lozenge Q}} \ar[r]^{\id} & TM\lozenge Q\\
T(M\lozenge Q)\ar[rr]^{T(\eta_M\lozenge Q)}&&
T(TM\lozenge Q)\ar[ru]_{U(\epsilon_{FM\tilde\lozenge Q})}, &
}$$
where the right-hand square is commutative by naturality 
of $\eta$ and the triangle represents the
adjunction identity $U(\epsilon_Z)\circ\eta_{UZ} = \id_Z$
for $Z = FM\tilde\lozenge Q$.

Axiom (D1) now follows by the commutativity of 
$$\xymatrix{TT(M\lozenge Q)
\ar[r]_{TT(\eta_M\lozenge Q)}\ar[dd]^{\mu_{M\lozenge Q}}&
TT(TM\lozenge Q)\ar[r]_{T(\widehat{\mu_M}^Q)}\ar[dd]^{\mu_{TM\lozenge Q}}&
T(TM\lozenge Q)\ar[r]_{T(\eta_{TM}\lozenge Q)}\ar[d]_{\id}&
T(TTM\lozenge Q)\ar[r]_{\widehat{\mu_{TM}}^Q}\ar[ld]^{T(\mu_M\lozenge Q)}
& TTM\lozenge Q \ar[dd]_{\mu_M\lozenge Q}\\
&&T(TM\lozenge Q)\ar[rrd]^{\widehat{\mu_M}^Q}
&&\\
T(M\lozenge Q)\ar[r]^{T(\eta_M\lozenge Q)}
& T(TM\lozenge Q)\ar[rrr]^{\widehat{\mu_M}^Q}&&&TM\lozenge Q
}$$
Here the right-hand square is clearly just following 
from a naturality of $\mu$, and the little triangle in the middle
from $\mu_M \circ \eta_{TM} = \id_M$. The middle pentagon
follows after translating all 4 nontrivial maps in terms of 
counit $\epsilon$ of the adjunction. Namely, using
$TM\lozenge Q = U(FM\tilde\lozenge Q)$, $\mu_X = U(\epsilon_{UX})$
and $\widehat{\mu_X}^Q = U(\epsilon_{FX\tilde\lozenge Q})$,
we can write the pentagon (after erasing $\id$-leg) as a square
$$\xymatrix{
UFUFU(FM\tilde\lozenge Q)\ar[d]^{U(\epsilon_{U(FM\tilde\lozenge Q)})}
\ar[rr]^{UFU(\epsilon_{FM\tilde\lozenge Q})} &&
UFU(FM\tilde\lozenge Q)\ar[d]^{U(\epsilon_{FM\tilde\lozenge Q})} \\
UFU(FM\tilde\lozenge Q)\ar[rr]^{U(\epsilon_{FM\tilde\lozenge Q})} 
&& U(FM\tilde\lozenge Q)
}$$
which is commutative by the naturality of $\epsilon$ 
(without $\tilde\lozenge Q$
everywhere this is nothing but the associativity diagram for $\mu$).
It is similar with the right-hand most rectangle which follows from
$$\xymatrix{
FU(FUFM\tilde\lozenge Q) \ar[r]^{\epsilon_{FUFM\tilde\lozenge Q}}
\ar[d]_{FU(\mu_M\tilde\lozenge Q)} & FUFM\tilde\lozenge Q 
\ar[d]^{\mu_M\tilde\lozenge Q} \\
FU(FM\tilde\lozenge Q)\ar[r]^{\epsilon_{FM\tilde\lozenge Q}} &
FM\tilde\lozenge Q,
}$$
after applying functor $U$ to every piece of the diagram and 
straightforward renamings; in this diagram $\mu_M$ is understood as
a map $\epsilon_{FM} : FUFM\to FM$  rather than $U(\epsilon_{FM})$.

To show (D2), one inspects in a similar manner this diagram:
$$\xymatrix{
T(M\lozenge(Q\otimes Q'))\ar[r]^{T(\eta_M\lozenge (Q\otimes Q'))}
\ar[d]^{T\Psi_M^{Q,Q'}} &
T(TM\lozenge(Q\otimes Q'))
\ar[rr]^{U(\epsilon_{FM\tilde\lozenge (Q\otimes Q')})}
\ar[d]^{T\Psi_{TM}^{Q,Q'}}&& 
TM \lozenge (Q\otimes Q')\ar[d]_{\Psi_{TM}^{Q,Q'}}\\
T((M\lozenge Q)\lozenge Q')\ar@/^/[r]^{T((\eta_M\lozenge Q)\lozenge Q')}
\ar[d]^{T(\eta_{M\lozenge Q}\lozenge Q')} &
T((TM\lozenge Q)\lozenge Q')\ar@/^/[rd]^{\id} 
\ar[rr]^{U(\epsilon_{(FM\tilde\lozenge Q)\tilde\lozenge Q'})}
\ar[d]^{T(\eta_{TM\lozenge Q}\lozenge Q')} &&
(TM\lozenge Q)\lozenge Q'\\
T(T(M\lozenge Q)\lozenge Q')\ar[r]^{T(T(\eta_M\lozenge Q)\lozenge Q')}
\ar[d]^{U(\epsilon_{FT(M\lozenge Q)\tilde\lozenge Q'})} &
T(T(TM\lozenge Q)\lozenge Q')
\ar[r]^{T(U(\epsilon_{FM\tilde\lozenge Q})\lozenge Q')}
\ar[rd]_{U(\epsilon_{F(TM\lozenge Q)\lozenge_l Q'})}
& T((TM\lozenge Q)\lozenge Q')
\ar@/_/[ru]^{U(\epsilon_{(FM\tilde\lozenge Q)\tilde\lozenge Q'})}& \\
T(M\lozenge Q)\lozenge Q' \ar[rr]^{T(\eta_M\lozenge Q)\lozenge Q'}&&
T(TM\lozenge Q)
\lozenge Q'.\ar[ruu]_{U(\epsilon_{FM\tilde\lozenge Q})\lozenge Q'}&
}$$

Finally, to show (D4), it is enough to consider the diagram
$$\xymatrix{ 
TM\ar@/^/[rrrr]^{\id} \ar[rr]_{T(\eta_M)}\ar[d]_{T(u_M)}&&
TTM\ar[rr]_{U(\epsilon_{FM})}\ar[d]^{T(u_{TM})}&&
TM\ar[d]^{u_{TM}}\\
T(M\lozenge 1)\ar@/_/[rrrr]_{l^Q_M}
\ar[rr]^{T(\eta_M\lozenge 1)}&&
T(TM\lozenge 1)\ar[rr]^{U(\epsilon_{FM\tilde\lozenge 1})}&&
TM\lozenge 1.
}$$

Finally, the fact that the correspondences $\tilde\lozenge \mapsto l$
and $l \mapsto \lozenge_l$ are inverses at one side has been shown above
(see (*)); and the other composition is left to the reader. 
Then the proof of the theorem is finished.

\vskip .15in
\ppt Given a commutative unital ring $\genfd$ and
a $\genfd$-coalgebra $C$ with coproduct 
$\Delta = \Delta^C : c\mapsto \sum c_{(1)}\otimes_\genfd c_{(2)}$ 
(Sweedler notation,~\cite{Majid})
we can consider the category
${\cal C} = {\cal C}^C$ of right $C$-comodules. 
 
Suppose now $B = C$ is a bialgebra. Then the category ${\cal C}^B$
is monoidal. Namely, if the coproduct on object $Q$ in ${\cal C}$ is given, 
in Sweedler-like notation, by $x \mapsto \sum q_{(0)}\otimes q_{(1)}$
then the tensor product in ${\cal C}$ is the tensor product 
$Q\otimes_\genfd Q'$ of $\genfd$-modules 
with the coaction $q\otimes q' \mapsto
\sum q_{(0)}\otimes_\genfd 
q'_{(0)}\otimes_\genfd q_{(1)}q'_{(1)}$. 
Bialgebra $B$ is a {\it monoid}
in ${\cal C}^B$, hence it induces a monad in any category on
which ${\cal C}^B$ acts. 
Category ${\cal C}^B$ acts on the category 
${\cal V}:= \genfd-{\rm Mod}$ of $\genfd$-modules
via the ``forgeftul'' right action, namely
$V\lozenge Q$ is $V\otimes_\genfd Q$ as a $\genfd$-module.
In particular, $B$ induces a monoid with the underlying
endofunctor $V \mapsto V\otimes_\genfd B$ in ${\cal V}$ as well.

A left $B$-module algebra $A = (A,\mu^A,\eta^A)$ 
is an associative unital algebra with
a left $B$-action $\triangleright_A$ for which 
$b\triangleright_A (a a') = \sum (b_{(1)}\triangleright_A a)
(b_{(1)}\triangleright_A a)$. Category $\tilde{\cal V} = A-{\rm Mod}$ 
of left $A$-modules is, in general, 
not monoidal (as $A$ is not commutative),
but it is equipped with the forgetful functor 
$F : \tilde{\cal V}\to {\cal V}$. 

Algebra $A$ induces a natural monad ${\bf T}= (T,\mu,\eta)$ 
in ${\cal V}$ given by 
$T(V) = A\otimes_\genfd V$ with $\mu_V = \mu^A \otimes_\genfd V$.
Clearly, ${\bf T}$-modules are the same thing as the left $A$-modules,
as $\nu : TM = A\otimes M\to M$ is a ${\bf T}$-action 
iff it is a left $A$-action.
\vskip .02in

We may consider the problem of lifting the forgetful action of
 ${\cal C}^B$ on ${\cal V}$ to 
${\bf T}-{\rm Mod} \equiv({\cal C}^B)^{\bf T}$.
There is a natural distributive law $l$ 
from $\lozenge$ to $\bf T$ compatible with $F$.
Namely given a module $M$ with action
$\triangleright_M : A\otimes M\to M$,
$l^Q_M : T(M\lozenge Q) = A\otimes (M\lozenge Q)\to
(A\otimes M)\lozenge Q) = TM\lozenge Q$ is given by
$a\otimes (m\otimes q) \mapsto 
\sum (q_{(1)}\triangleright_A a \otimes m)\otimes q_{(0)}$.
Hence the induced action is
$(M,\nu)\lozenge_l Q = (M\otimes Q,\nu^Q)$ where 
$a\triangleright_{TM} (m\otimes q) \equiv \nu^Q(a\otimes (m\otimes q)) = 
\sum (q_{(1)}\triangleright_A a)\triangleright_M m)\otimes q_{(0)}$.

Consider the associated monoidal functor 
${\cal L} : {\cal C}\to {\rm End}\,\tilde{\cal V}$
for the lifted action. 
As the bialgebra itself is a monoid in ${\cal C}^B$, its image
is also a monoid in ${\rm End}\,\tilde{\cal V}$, 
that is a monad in $\tilde{\cal V}$. 
Call this new monad ${\cal L}_B$.
Then ${\cal L}_B (M,\nu) = (M\otimes B,\nu^B)$. 
Its multiplication 
$\mu' = \mu^{{\cal L}^B} : M\otimes B\otimes B\to M\otimes B$
is simply $M\otimes \mu^B$. 
We leave as an exercise for the reader to analyse that 
the Eilenberg-Moore category $\tilde{\cal V}^{{\cal L}^B}$
of this monad is equivalent to a certain category of $\genfd$-modules
equipped with two additional actions (the left action of $A$
inherited from $\cal V$ and a new right action of $B$)
with the compatibility condition 
$a \triangleright (n \blacktriangleleft h) = 
[(h_{(2)} \triangleright_A a) \triangleright n] \blacktriangleleft h_{(1)}$.
$\tilde{\cal V}$ is not monoidal, and even if it is 
(if $A$ is commutative, say),
monad ${\cal L}_B$ is {\it not}
induced by tensoring by an algebra in $\tilde{\cal V}$. 
Instead, $B$ is an algebra in another monoidal category -- namely in
${\cal C}^B$; but also in the underlying category ${\cal V}$.

This simple example has a number of useful generalizations and applications
which will be addressed elsewhere. 

\ppt There are other flavours
of compatibility between tensor structures and monads. 
Instead of ${\cal C}$-categories in the sense of actions, one can consider
the different notion of categories enriched over ${\cal C}$ (i.e. the
hom-sets consist of $\cal C$-objects, and carry tensor products 
with some natural properties), those are also called
${\cal C}$-categories; monads in such categories 
were explored before~(\cite{Bunge:relfunct,Wolff:locmon}). 
One can also consider the monads 
in ${\cal C}$-itself, cf.~\cite{McCrudden,Moerd:tenmon}.
There are also various dual 
(left-hand, opmonoidal, comonadic etc.) versions. 

\ppta Such a ``dual'' version which may be less obvious is as follows.
Given ${\cal C}$-category ${\cal M}$  
in our sense (i.e. with action $\lozenge$), 
one considers the category ${\cal M}^{\lozenge}$
of ${\cal C}$-modules {\it in} 
${\cal M}$, i.e. objects $M$ in ${\cal M}$ with families of maps
$\nu_Q : M\lozenge Q\to M$ satisfying a list of required identities
making those famililes ``actions''. 
Such modules also make a category 
in a natural way and one can consider liftings of a
monad ${\bf T}$ in ${\cal M}$ to that category. 
The liftings will then be given by another flavour of distributive laws. 

\vskip .2in
\footnotesize{{\bf Acknowledgements}. The main result has been obtained
during author's visit to IH\'ES whose hospitality in Spring 2004 resulted in
a beautiful and stimulating scientific experience. 
I also thank the organizers
of ``NOG III, Quantum geometry'' conference at Mittag Leffler, particularly
Prof. Laudal and Prof. Kontsevich for the invitation to present 
a talk where a wider context pertaining to this 
and related results has been exhibited. I thank Prof. Brzezi\'{n}ski for
finding some expositional errors.
}

\end{document}